\documentclass{article}[10pt]

\usepackage{latexsym}
\usepackage{amssymb}
\usepackage{amsmath,amsfonts}
\usepackage{tikz}
\usepackage{graphicx}
\usepackage{multicol}
\usepackage{multirow}
\usepackage{subfigure}
\usepackage{color}
\usepackage{graphicx,caption}
\usepackage{latexsym,booktabs}
\usepackage{amsthm}
\usepackage{mathrsfs}
\usepackage{caption}

\tikzstyle{vertex}=[circle, draw, inner sep=0pt, minimum size=6pt] 
\newcommand{\vertex}{\node[vertex]}
\usetikzlibrary{decorations.markings}

\def \Q {{\mathbb Q}}
\def \N {{\mathbb N}}
\def \Z {{\mathbb Z}}
\def \R {{\mathbb R}}
\def \F {{\mathbb F}}
\def \E {{\mathbb E}}
\def \L {{\mathbb L}}
\def \C {{\mathbb C}}
\def \K {{\mathbb K}}
\def \u {{\bf u}}
\def \0 {{\mathbf 0}}
\def \B {{\mathcal B}}

\def \la {{\lambda}}

\newtheorem{theorem}{Theorem}[section]

\theoremstyle{definition}
\newtheorem{defn}{Definition}[section]

\newtheorem{ex}{Example}[]
\newtheorem{pro}[theorem]{Proposition}
\theoremstyle{remark}

\newtheorem*{note}{Note}

\def\[{[\hskip-1pt [}
\def\]{]\hskip-1pt ]}
\def\N{{\mathbb N}}

\def\++{\boxplus}

\begin{document}

\baselineskip 15pt

\title{K-Primitivity : A Literature Survey}
\date{}
\author {Monimala Nej\\
Department of 
Mathematics \\ Adamas
University, India-700126.
  }
\maketitle
\begin{abstract}
A nonnegative matrix A is said to be {\it primitive} if there exists a positive integer $m$ such that entries in $A^m$ are positive and smallest such $m$ is called the {\it exponent} of $A.$ Primitive matrices are useful in the study of finite Markov chains theory. In 1998, in the context of finite Markov chains, Ettore Fornasini and Maria Elena Valcher \cite{Forn:Val} extended the notion of primitivity for a nonnegative matrix pair $(A,B)$ by considering a positive discrete homogeneous two-dimensional $(2D)$ state model. Further generalization to this notion of primitivity for $k$-tuple $(A_1,A_2,\ldots,A_k)$ of nonnegative matrices $A_1, A_2, \ldots, A_k$ is quite natural and known as {\it $k$-primitivity}. In this paper we present various results on $k$-primitivity given by different researchers from time to time.
\end{abstract}
{\bf{AMS(2020)}}: 05C50, 05C38, 15B99.
\section{Introduction}\label{sec:1}
Let $A\in M_n(\R),$ where $M_n(\R)$ is the set of $n\times n$ matrices with real entries. 
If  all  the  entries of $A$  are positive, then $A$  is called a {\it 
positive} 
matrix. We denote a positive matrix as $A>0$. Similarly if each entry of a real 
matrix $A$ is nonnegative, then we denote it as $A\ge 0$ and  is called {\it 
nonnegative} matrix. Nonnegative matrices and their eigenvalues have applications in  several areas of mathematics. One such important area in our context is probability theory,  in particular, finite Markov chains. In Markov chain theory, regular chain is one whose transition matrix is a primitive matrix. But there, a primitive matrix usually called as a regular matrix\index{regular matrix}, refer to Peter Perkins \cite{Per} and G.L. Thompson \cite{Thhhom}.\\
Let us consider a Markov chain
$$x_k = Ax_{k-1} \; (k \geq 1),$$ where A is a nonnegative $n \times n$ matrix and the initial value $x_0$ is a nonzero, nonnegative vector. It has proven useful to distinguish Markov chains for which $x_k > 0$ for all sufficiently large $k.$ By considering the exponential growth rate of the matrix power $A^k$, it can be found that for each nonzero, nonnegative initial condition $x_0$ there is a $K$ such that $x_k > 0$ for all $k \geq K$ if and only if there exists a positive integer $l$ such that $A^l > 0.$ Further, the exponential growth rate of the matrix power $A^k$ as $k \to \infty$ is controlled by the eigenvalue of $A$ with the largest absolute value. Thus, primitivity can be viewed as a property of a nonnegative matrix in spectral theory and on other hand it can be viewed as an important property of a nonnegative matrix in Markov chain theory.

We now discuss the definition of {\it primitive matrix} which arises from the characterization of the eigenvalues of the matrix. Then we see an alternative definition for the same which can be found in Markov chain theory and the key point for defining the $k$-primitivity of $k$-tuple $(A_1,A_2,\ldots,A_k)$ of nonnegative matrices. A $k$-tuple of nonnegative matrices which is $k$-primitive also carry various spectral properties similar to a primitive matrix. Let us see a few necessary definitions, results in  the development of the nonnegative matrix theory, specially on eigenvalues and eigenvectors.

In 1907, Oskar Perron proved that a real square matrix with positive entries has a unique largest real eigenvalue and that the corresponding eigenvector can be chosen to have strictly positive components.

Let $A\in M_n(\R)$ and $\la_1,\la_2,\ldots,\la_n$ be eigenvalues of $A$ with repetitions, then {\it spectral radius\index{spectral radius} of $A$} denoted $\rho(A)$ and is defined as 
$$\rho(A)=\max\{|\la_1|,|\la_2|,\ldots,|\la_n|\}.$$ Note that $\rho(A)$ need not be an eigenvalue of $A.$ For example, if $A=\begin{bmatrix}
              -1 & 1\\
              0&-1
                                                                                                                             \end{bmatrix},$ then $\rho(A)=1$ and  is not an eigenvalue of $A.$ Hence for a matrix $A > 0$ we have the following result.

                                                                                                                 \begin{pro}[O. Perron \cite{Perr}, 1907]\label{thm:perron} Let $A$ be a positive matrix and $r=\rho(A).$ Then
\begin{enumerate}
\item $r$ is an eigenvalue of $A$ called the {\it Perron root}\index{Perron root} of $A.$ If $\la$ is any other eigenvalue of $A$, then $|\la|<r.$
\item  algebraic and geometric multiplicity of $r$ is one.
\item there is a positive eigenvector corresponding to $r.$ Moreover, there are no other nonnegative eigenvectors of $A$ other than the vectors in the eigenspace associated to $r.$
\end{enumerate}
\end{pro}
Later, in 1912, George Frobenius extended the results in Proposition \ref{thm:perron} to the class of irreducible matrices, a particular class of nonnegative matrices. He described the properties of the leading eigenvalues and the properties of the corresponding eigenvectors when $A \geq 0$. To proceed further, let us see the following definition. 
\begin{defn}[H. Minc \cite{Minc}]
\begin{enumerate}
\item
A matrix $A$ is said to be {\it permutationally similar} to a matrix $B$, denoted $A \backsimeq_PB$, if there exists a permutation matrix $P$ satisfying $B=PAP^t$, where $P^t$ denotes the transpose of $P.$
\item
A nonnegative $n$-square matrix $A$, $n \geq 2$, is called {\it reducible(decomposable)} if there exists a permutation matrix $P$ such that $$PAP^t= \begin{bmatrix}
               B& C\\
               0 &D 
              \end{bmatrix},$$ where $B$ and $D$ are square submatrices. Otherwise $A$ is {\it irreducible(indecomposable)}.
              
\end{enumerate}
\end{defn}
\begin{ex}\label{exam:thesis:irre}
\begin{enumerate}
\item
The identity matrix $I_n, n \geq 2$ is reducible. 
\item If $A=\begin{bmatrix}
0&1&0\\
0&0&1\\
0&0&0
\end{bmatrix},$ then $A$ is reducible as $A=\begin{bmatrix}
               B& C\\
               0 &D 
              \end{bmatrix},$ where\\ $B=[0], C=[1\;\; 0]$ and $D=\begin{bmatrix}
0&1\\
0&0
\end{bmatrix}.$ 
\item Any positive matrix is irreducible. A simple nontrivial irreducible matrix is $\begin{bmatrix}
0&1\\
1&0
\end{bmatrix}.$ 
\end{enumerate}
\end{ex} 

If $A$ is a nonnegative matrix of the form $\begin{bmatrix}
               B& C\\
               0 &D 
              \end{bmatrix},$ where $B$ and $D$ are square matrices, then $A^n$ is also of the same form for every $n.$ Hence if $A$ is reducible, then $A^n$ is also reducible for every $n\in \N.$ Thus a necessary condition for   a nonnegative matrix  to be primitive is that it is irreducible. The following result can be used to check whether a given matrix is irreducible or not.
\begin{pro} [H. Minc \cite{Minc}] \label{pro:irr} Let $A$ be a nonnegative matrix of order $n.$ Then $A$ is irreducible if and only if $(I_n+A)^{n-1}>0.$
\end{pro}
In 1912, G. Frobenius established the following results in Proposition \ref{thm:perron-frob} which have important applications to probability theory, specially to Markov chains.
\begin{pro}[G. Frobenius \cite{Frob:2}, 1912]\label{thm:perron-frob}
Suppose $A$ is an irreducible matrix. Then
\begin{enumerate}
\item
$\rho(A)$ is an eigenvalue of $A.$ For any other eigenvalue $\lambda$ of $A$, $|\lambda| \leq \rho(A).$
\item
the algebraic multiplicity and geometric multiplicity of $\rho(A)$ is $1.$
\item
there is a positive eigenvector corresponding to $\lambda.$ Moreover there are no other nonnegative eigenvectors of $A$ other than the vectors in eigenspaces associated to $\rho(A).$
\end{enumerate}
\end{pro}

Now the results of Proposition \ref{thm:perron-frob} lead for the following definition. And Theorem \ref{prim:expo} can be established by using this definition.
\begin{defn}[see, H. Minc \cite{Minc}]\label{defn:1}
Let $A$ be an $n \times n$ irreducible matrix with $r=\rho(A)$. Suppose that $A$ has exactly $h$ eigenvalues of modulus $r.$ The number $h$ is called the {\it index of imprimitivity} or simply the {\it index} of $A$. If $h=1$, then the matrix $A$ is said to be {\it primitive} otherwise it is {\it imprimitive} (or cyclic, in some authors' nomenclature).
\end{defn}
Thus, due to G. Frobenius \cite{Frob:2}, an alternative definition of primitive matrix can be found in Theorem \ref{prim:expo}.
\begin{theorem}\label{prim:expo}
A necessary and sufficient condition for a nonnegative matrix $A$ to be primitive is that $A^m >0$ for some positive integer $m.$
\end{theorem}
Let $A$ be an $n \times n$ nonnegative matrix which is primitive. Then $\min \{m: A^m >0\}$ is called the {\it primitive exponent} or simply the {\it exponent} of $A$ and denoted $exp(A)$ or $\gamma(A).$ Sometimes it is also called as {\it index of primitivity}, see B.R. Heap and M.S. Lynn \cite{Heap:Lynn}. In 1954, I.N. Herstein \cite{Herstein} gave a simple proof of Theorem \ref{prim:expo} which is algebraic in nature and avoids the use of the convergence of powers of a matrix.

Research focused on primitive exponent ever since 1950, when Wielandt published his paper \cite{Wie}. Since then various results are established in this context. A Brief survey on the same can be found in J. Shen \cite{Shen:9}, B. Liu and Hong-Jian Lai \cite{Liu:Lai} (p. 107-114). In 1997 E. Fornasini and M.E. Valcher \cite{Foorrn:Vaaaalll} generalized the notions of irreducibility, primitivity of a nonegative matrix to $k$-tuple of nonnegative matrices of the same order. In this paper we discuss various results on the notion of {\it k-primitivity}. Before further continuation, we will introduce some basic concepts on graph theory. The theories of graphs play an important role to esablish various results on primitivity as well as on $k$-primitivity of nonnegative matrices. 
\section{Graph Theory}
It is known that every matrix with $0$ and $1$ entries corresponds to a  directed graph and vice versa. After introducing a few definitions and notations of graph theory, we will see  graph theoretic interpretations of irreducible matrices and primitive matrices. We denote $D=(V,E)$ as a {\it digraph}\index{digraph} (directed graph) with vertex set\index{vertex set} $V=V(D)$ and the edge set\index{edge set} $E=E(D)$ and order $n=|V|.$  A digraph we mean no multiple edges but loops are allowed. A $u \to v$ {\it walk}\index{walk} in $D$ is a sequence of vertices $u, u_1, u_2, \ldots, u_l=v$ and a sequence of edges $(u,u_1), (u_1,u_2), \ldots, (u_{l-1},v)$ where vertices and edges may be repeated. A {\it cycle} or {\it closed walk}\index{closed walk} is a $u \to v$ walk where $u=v.$ A {\it path}\index{path} is a walk with distinct vertices. An {\it elementary cycle}\index{elementary cycle} or {\it circuit}\index{circuit} is a closed $u \to v$ walk with distinct vertices except for $u=v.$ The length of a walk $W$ is the number of edges in $W.$ A circuit of length $1$ is nothing but a loop, and conversely.

Let $D$ be a digraph on $n$ vertices and let us fix a
labeling of the
vertices of  $D$. Then, the {\it adjacency
matrix}\index{adjacency matrix}
of $D$, denoted $A(D) = (a_{i,j})$
(or $A$), is an $n \times n$ matrix with $a_{i,j} = 1$, if
the vertex  $i$ is
adjacent to the vertex $j$ and $0$, otherwise.
Let $D$ be a directed graph and let $B$ be its adjacency matrix. Then,  $D$ is 
said to be associated with a nonnegative matrix $C$ if
$C_{i,j}\ne 0\Leftrightarrow B_{i,j}\ne 0.$  We denoted the digraph $D$ associated 
with $C$ as $D(C).$ In our context, we call $D(C)$ as the {\it digraph}\index{digraph} of $C.$ If
$C$ is a symmetric matrix,
then $D(C)$ is an undirected graph, or simply a {\it graph}\index{graph}.

Note that a single digraph can be associated to more than one matrix. For example, a 
complete graph\index{complete graph} with $|V|=n$ and a loop at each vertex is associated to every 
positive matrix, {\it i.e.}, adjacency matrix of $D(C)$ is $J_n$ whenever $C>0.$

\begin{ex} \label{nolabel}
The digraphs of $A$ and $B$ are same, where 
\end{ex}

\begin{tabular}{lll}
$A =
\begin{bmatrix}
0 & 1 & 0 & 0\\
0 & 0 & 1 & 0\\
1 & 0 & 0 & 1\\
1 & 0 & 0 & 1
\end{bmatrix},$  & $B=\begin{bmatrix}
0 & 1 & 0 & 0\\
0 & 0 & 2 & 0\\
7 & 0 & 0 & 3\\
4 & 0 & 0 & 5
\end{bmatrix} \;\;and $\\
\end{tabular}

\begin{figure}[htb]
\centering
 \begin{tikzpicture}
 [scale=1]
       \vertex (c1) at (90:1) [label=above:$1$]{};
       \vertex (c2) at (180:1) [label=above:$2$]{};
       \vertex (c3) at (270:1) [label=left:$3$]{};
	\vertex (c4) at (360:1) [label=below:$4$]{}edge [in=80,out=300,loop] ();

	\path[->]
	        (c1) edge (c2)
		(c2) edge (c3)
		(c3) edge (c4)
		(c3) edge (c1)
		(c4) edge (c1)
	
		;
		
\end{tikzpicture} 
\caption*{$D(A)$ or $D(B).$}
\end{figure}
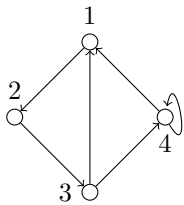

If $A$ is the adjacency matrix of a digraph $D$ then it is known that the $ij$-th entry of $A^k$ denotes the number of walks of length $k$ from the vertex $i$ to $j.$ 

Let $B$ be an $n \times n$ irreducible matrix, then from the Proposition \ref{pro:irr} for every $i$, $j,$ $1 \leq i,j \leq n,$ there exists $k$ such that $ij$-th entry of $B^k$ is positive. That is, there is a path  between any two vertices in $D(B).$ Digraphs with this property are called {\it strongly connected}. Thus a digraph is strongly connected if and only if its adjacency matrix is irreducible. Further, if $B$ is a primitive matrix, then there exists $m \in \N$ such that $B^m>0.$ Equivalently, there is a walk of length $m$ between any two vertices in $D(B)$. With these observations, a digraph is called irreducible (primitive)\index{primitive digraph} if its adjacency matrix is irreducible (primitive). Hence the graph theory and the matrix theory can be used interchangeably on the study of primitive matrices and exponents. The following result is useful in this direction.

\begin{pro} [V. Romanovsky \cite{Romanovskyyy}]
A digraph $D$ is primitive if and only if $D$ is strongly connected and the g.c.d. of the lengths of all circuits in $D$ is $1$. 
\end{pro}

\begin{ex}
It is easy to see that the exponent of the adjacency matrix $A$ of the star graph 
$K_{1,n}$ with a loop at the vertex of degree $n$ is $2$. The graph $K_{1,5}$ with a 
loop at the
vertex $1$ of degree $5$ is drawn below. There are $6$ walks of length $2$ from the vertex $1$ to itself. And the number 
of walks of length $2$ for the remaining combinations is $1$. 
\end{ex}
\begin{figure}[htb]
\centering
 \begin{tikzpicture}
	\vertex (c1) at (89:2) [label=above:$6$]{};
	\vertex (c2) at (161:2) [label=above:$5$]{};
	\vertex (c3) at (233:2) [label=left:$4$]{};
	\vertex (c4) at (305:2) [label=below:$3$]{};
	\vertex (c5) at (17:2) [label=right:$2$]{};
	\vertex (c6) at (0,0) [label=below:$1$] {};
	\path
                (c6) edge (c1)	
		(c6) edge (c2)
		(c6) edge (c3)
		(c6) edge (c4)
		(c6) edge (c5)
		edge [in=140,out=60,loop] ()

		;
		
\end{tikzpicture}
\caption{The star graph $K_{1,5}$ with a loop at the vertex $1$}\label{star:loop:new}
\end{figure}
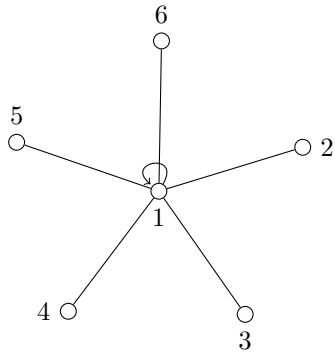
\begin{tabular}{lll}
$A= \begin{bmatrix}
1 &1 &1&1&1&1\\
1&0&0&0&0&0\\
1&0&0&0&0&0\\
1&0&0&0&0&0\\
1&0&0&0&0&0\\
1&0&0&0&0&0
\end{bmatrix}$
& $A^2= \begin{bmatrix}
6 &1 &1&1&1&1\\
1&1&1&1&1&1\\
1&1&1&1&1&1\\
1&1&1&1&1&1\\
1&1&1&1&1&1\\
1&1&1&1&1&1
\end{bmatrix}.$\\
\end{tabular}\\
\  \\
\ \\

Let $\B_n$ be the set of all matrices of order $n$ with entries $0,$ $1$ and $1+1=1$ in the addition of the entries in these matrices. From above discussion, working with nonnegative matrices in the context of primitivity, it is sufficient to work with matrices from $\B_n.$  Clearly there is a one-to-one correspondence between matrices in $\B_n$ and all labeled digraphs with $n$ vertices. Hence, for rest of this paper, graph theory and matrix theory are used interchangeably for convenience. Occasionally rest of the definitions and notations will be given at the appropriate place.

\section{Concept of $k$-primitivity and relevant results}
Primitive matrices plays an important role in the study of finite Markov chains theory. In finite Markov chains, different kinds of discrete models are studied. One such simple homogeneous discrete model is $$x(t+1) = A x(t), \hspace{5mm} t = 0, 1, \ldots \hspace{50mm}(1)$$ where $x(k)$ represent the $k$-th state of the model and $x(0)$ is the given initial condition, $A$ is a nonnegative matrix associated with the model. It is well-known that $A$ is said to be primitive if and only if for every nonzero, nonnegative initial condition $x(0),$ the state of the said model becomes positive in a finite number of steps. Later, this notion of primitivity has been extended to a nonnegative matrix pair $(A,B)$ by considering a discrete homogeneous $2D$-system, namely, 
$$x(h+1, k+1) = Ax(h, k+1) + Bx(h+1, k), \hspace{3mm} h, k \in \Z, \; h + k \geq 0, \hspace{8mm}(2)$$
where $x(h,k)$ denotes the local states of the system and $x(h, -h) \; (h \in \Z)$ are given initial conditions which are nonnegative $n \times 1$ matrices and $A$, $B$ are $n \times n$ nonnegative matrices. Now the initial values $x(h,-h)$ are said to be admissible provided there exists an integer $N$ such that $\sum\limits_{i=0}^N x(h+i,-h-i)$ is nonzero and nonnegative for all $h.$

Various properties of pairs of $n \times n$ nonnegative matrices are useful in the context of certain dynamical systems that generalize Markov chains. In 1998 E. Fornasini and M.E. Valcher  \cite{Forn:Val} investigated a dynamical $2D$ system via a directed graph and extended several important spectral properties of nonnegative matrices to pairs of nonnegative matrices. One such important spectral property is primitivity which is defined for a nonnegative matrix pair $(A,B)$ as follows.

\begin{defn}(\cite{Forn:Val},1998)
A pair of nonnegative matrices $(A,B)$ is primitive if for every admissible sequence of initial conditions the local state $x(h,k)$ becomes strictly positive when $h+k$ is sufficiently large. 
\end{defn}
Here we would like to mention that primitivity of circulant matrices plays an important role to the primitivity of the pair $(A,B)$ whenever the initial conditions of the associated $2D$-system are periodic. Similarity carried out for the non-periodic initial conditions also. Now we see an alternative definition for the primitivity of the pair $(A,B)$ using Hurwitz product. This definition also can be found in \cite{Forn:Val}.
\begin{defn} (see, for instance, \cite{Bea:Mou})
Let $A_1 , A_2 , \ldots , A_k$ are nonnegative square matrices of order $n$ and let $m_1, \ldots m_k$ be positive integers. The sum of all possible products containing $m_i$ $A_i$'s for all $i \in \{1, \ldots ,k\}$ is called the {\it $(m_1, \ldots ,m_k)$-Hurwitz product} of $A_1, \ldots ,A_k.$ It is denoted by $(A_1, \ldots ,A_k )^{(m_1, \ldots ,m_k)}.$ 
\end{defn}
For example, $(2,2)$-Hurwitz product of $(A,B)$ is
$$(A,B)^{(2,2)} = AABB+BBAA+BABA+ABAB+ABBA+BAAB.$$ Now $(A_1, \ldots ,A_k )^{(m_1, \ldots ,m_k)} > 0$ will indicate that $(m_1, \ldots ,m_k)$-Hurwitz product of $(A_1 , A_2 , \ldots , A_k)$ is positive (entrywise).
\begin{defn}
A pair of nonnegative matrices $(A,B)$ is said to primitive if there exists nonnegative integers $h$ and $k$ such that $h+k > 0$ and $(A,B)^{(h,k)} > 0.$ And the exponent of $(A, B)$ is defined to be the minimum of $h+k$ over all ordered pair $(h,k)$ of nonnegative integers for which $(A,B)^{(h,k)} > 0.$ This is called the {\it $2$-exponent} of $2$-primitive matrix pair $(A,B)$.
\end{defn}
It is quite natural to extend this notion of primitivity to any $k$-tuple ($k \in \N$) $\mathscr{A}=(A_1,A_2,\ldots,A_k)$ of nonzero and nonnegative matrices of same order. In this case $\mathscr{A}$ is said to be primitive or very often {\it $k$-primitive.} And the exponent of $\mathscr{A}$ is defined to be the minimum of $i_1+i_2+ \ldots + i_k$ over all nonnegative integer $k$-tuples $(i_1, i_2, \ldots, i_k)$ for which $\mathscr{A}^{(i_1,i_2 ,\ldots,i_k)}>0.$ This is called the {\it $k$-exponent} of $k$-primitive matrix tuple $\mathscr{A}$ and denoted by $exp(\mathscr{A})$.

In this context of primitivity, the graph is defined for $(A,B)$ which is called the multidigraph and denoted by $D^2.$ If the order of the matrix $A$ or $B$ is $n,$ then $D^2$ is a graph with $n$ vertices and arcs of two different kinds, namely $A$-arcs and $B$-arcs. There is an $A$-arc from vertex $i$ to vertex $j$ if and only if  $a_{ij} > 0,$ and a $B$-arc if and only if  $b_{ij} > 0.$ In general, consider $k$-tuple $\mathscr{A}=(A_1,A_2,\ldots,A_k)$ of nonzero and nonnegative matrices of same order $n$ and let $c_1, c_2, \ldots ,c_k$ be distinct colors. The {\it colored multidigraph} $D(\mathscr{A})$ of $\mathscr{A}$ is the multidigraph with vertices $1, 2, \ldots, n,$ and an arc of color $c_i$ from vertex $u$ to vertex $v$ if and only if the $(u,v)$-entry of $A_i$ is positive. Sometimes this digraph is called as {\it $k$-colored digraph} and denoted as $D^k$. A directed walk in $D(\mathscr{A})$ from $u$ to $v$ is a sequence $a_1, a_2 , \ldots, a_l$ of arcs (irrespective of color) of $D(\mathscr{A})$ such that the initial vertex of $a_1$ is $u,$ the terminal vertex of $a_l$ is $v$ and for $i = 1, 2, \ldots , l-1$ the terminal vertex of $a_i$ is the initial vertex of $a_{i+1}.$ With these terminologies, primitivity of $D^k$ has been defined as follows.

\begin{defn}
A $k$-colored digraph $D^k$ with distinct colors $c_1, c_2, \ldots ,c_k$ is said be {\it primitive} or {\it $k$-primitive} if there is a $k$-tuple of positive integers $(r_1, r_2 , \ldots , r_k)$ such that between any pair of vertices $u,v$ of $D,$ there is a walk from $u$ to $v$ in $D$ having exactly $r_i$ arcs of color $c_i$ for each $i,$ $i = 1, \ldots , k.$ 
\end{defn}

It is proved that a $k$-tuple $\mathscr{A}=(A_1,A_2,\ldots,A_k)$ of nonzero and nonnegative matrices of same order is primitive if and only $D^k$ is primitive. Hence $exp(\mathscr{A})=exp(D^k),$ where $exp(D^k)$ denotes the $k$-exponent of $D^k$ (see, \cite{YubinGaoYanlingShao}). In light of this correspondence one can freely move between two-colored digraphs and the associated matrix pairs. Thus primitivity and exponent of $\mathscr{A} = (A_1, A_2, \ldots ,A_k)$ depend only upon which entries of the $A_i$ are nonzero, and hence there is no loss of generality in restricting our study to $k$-tuples of $(0,1)$-matrices.  In 1997 E. Fornasini and M.E. Valcher \cite{Foorrn:Vaaaalll} introduce $2D$ digraph, namely a directed graph with two kinds of arcs. In this paper various spectral and combinatorial properties of a pair $(A,B)$ in a connection with associated $2D$ digraph have been studied and defined primitivity as a special case of irreducibility.  An irreducible matrix is said primitive if its {\it $2D$-imprimitivity index} $h^{(2)}(D)$ (see \cite{Foorrn:Vaaaalll}) is $1$ and characterized it in several alternative ways by resorting to the results derived in the same paper.


After E. Fornasini and M.E. Valcher, B.L. Shader, S. Suwilo \cite{Shad:Suw} were first to describe how certain $2D$-dynamical systems lead to an extension of the notions of primitivity and exponent to pairs of nonnegative matrices. They relate the exponent of $(A,B)$ to integers solutions to certain systems of linear diophantine equations and established that the exponent of a primitive pair $(A,B)$ of nonnegative $n \times n$ matrices is at most $\frac{(3n^3+2n^2-2n)}{2}.$  Further, they define a natural two-colored digraph, namely, {\it two-colored Wielandt digraph}, a two-colored digraph obtained by coloring the arcs of a Wielandt digraph. And showed that $2n^2-4n+1 \leq exp(D^{2}_W) \leq 2n^2-3n+1,$ where $D^{2}_W$ is a primitive, two-colored Wielandt digraph on $n \geq 3$ vertices with at least one red arc and one blue arc. Digraphs for the equality in either case of lower bound and upper bound can be found in S. Suwilo \cite{Suwiloooos}. In 2002, D.D. Olesky, B. Shader, P. van den Driessche \cite{Ole:Shad:Dri}, showed that for each positive integer $k,$ the maximum exponent of a primitive $k$-tuple of $n \times n$ nonnegative matrices is $\Theta(n^{k+1}).$

A new approach for characterizing $k$-primitive matrix families has been developed by V. Yu. Protasov \cite{Port}. In 2013, he proved that under some mild assumptions, a set of $k$ nonnegative matrices is either $k$-primitive or there  exists a nontrivial partition of the set of basis vectors in $\R^{n}$, on which these matrices act as commuting permutations. This gives a convenient classification of $k$-primitive families and a polynomial-time algorithm to recognize them. Similar classification of primitive family of nonnegative matrix tuples can be found in \cite{Prot:Voy}. A $k$-tuple $\mathscr{A} = (A_1, A_2, \ldots ,A_k)$ of nonnegative matrices of same order is said to be {\it primitive} if there is at least one strictly positive product (with repetitions permitted) of matrices from $\mathscr{A}$. This concept provides a
different generalization of primitivity of one matrix to matrix families. Clearly, primitivity
implies $k$-primitivity, but not vice versa. That is, primitivity is more generalized concept than $k$-primitivity for a $k$-tuple $\mathscr{A} = (A_1, A_2, \ldots ,A_k)$ of nonnegative matrices.


\subsection{$K$-Primitivity of a digraph $D$}
In continuation of research in this topic,  people obtained a uncolored digraph $D$ from a $k$-colored digraph $D^k.$ If $D^k$ is primitive, then $D$ is also primitive, and $exp(D) \leq exp(D^k).$ Suppose that $D$ is a simple directed graph on at least two vertices, possibly with loops. A {\it $k$-coloring of $D$} is a $k$-tuple of spanning subgraphs $(D_1 , \ldots , D_k )$ such that the subgraphs $D_1, \ldots ,D_k$ partition the arcs of $D$ into $k$ (nonempty) subsets. Each $D_i$ containing the arcs given color $i$ in the directed graph $D,$ see, for instance, \cite{LeeRoyy:Kkirrk}. Now a digraph $D$ is said to be $k$-primitive if there exists a $k$-coloring $(D_1 , \ldots , D_k )$ of $D$ such that there is a $k$-tuple of positive integers $(r_1, r_2 , \ldots , r_k)$ and between any pair of vertices $u,v$ of $D,$ there is a walk from $u$ to $v$ in $D$ having exactly $r_i$ arcs in $D_i$ for each $i,$ $i = 1, \ldots , k.$ Equivalently, that walk has exactly $r_i$ arcs of color $i$ for each $i = 1, \ldots , k.$ It is clear that if there is a $k$-coloring of $D$ that is $k$-primitive, then $D$ itself must be a primitive directed graph. On other hand, if $A_1, A_2 , \ldots , A_k$ be $(0,1)$-matrices such that their sum $A$ is also a $(0,1)$-matrix, that is, each nonzero entry of $A$ corresponds to exactly one nonzero entry in some $A_i,$ then it is known that $A$ is the disjoint sum of $A_1, A_2, \ldots , A_k.$ This decomposition of $A$ corresponds in a natural way to a $k$-coloring of the digraph $D.$ Now $A$ is said to be $k$-primitive if $k$-tuple $(A_1, A_2, \ldots , A_k)$ is $k$-primitive. Thus $A$ is $k$-primitive if and only if $D$ is so. 

It is known from \cite{Ole:Shad:Dri} that if digraph  $D$ is $k$-primitive, then necessarily
$D$ must have at least $k$ cycles. In 2003, LeRoy B. Beasley a , Steve Kirkland \cite{LeeRoyy:Kkirrk} showed that any primitive directed graph is a $2$-primitive and for each $k \geq 4,$ they provided an example of a primitive directed graph having exactly $k$ cycles, but not $k$-primitive. For $k=3,$ they provided one conjecture whose confirmation would be a key step in proving that every primitive directed graph with at least three cycles admits a $3$-coloring, that is $3$-primitive. This conjecture is not yet proved. In 2014, LeRoy B. Beasley and S. Mousley \cite{Bea:Mou} established conditions for primitive digraphs to be $k$-primitive which is basically an improvement of result that any primitive directed graph is a $2$-primitive. In the same paper they proved that $\binom{n-1}{2} \leq k_{max}(n) < \binom{n}{2}-\lceil \frac{n}{4} \rceil$, where $k_{max}(n)$ be the maximum $k$ for which there exists a $k$-coloring of some strongly connected $n$-tournament that is $k$-primitive. This result is an extention of results in \cite{Beassss:Neeealll}, where properties of 2-primitive tournament digraphs have been discussed. In 2009, Y. Gao, Y. Shao \cite{Gaooo:Shaaaoooo} showed that there are $k$-colorings that are $k$-primitive for two different class of digraphs, namely, the class of primitive symmetric digraph of order $n$ and for the class of primitive digraphs such that for each primitive digraph $D$ if there exist cycles $C_1, C_2 , \ldots , C_k,$ then $C_i$ and $C_j$ have no common arc for each pair $(i,j)$ with $1 \leq i <j \leq k.$

\subsection{Results on $2$-primitivity and $2$-exponent of $(A,B)$}
As of now, in the context of $k$-primitivity, most of the researchers focused on the study of $2$-primitivity and $2$-exponent. In this section we represent a few important results on $2$-primitivity and $2$-exponent for various family of $2$-tuples $\mathscr{A}=(A_1,A_2)$ of nonnegative matrices.

Let $D^2$ be a $2$-colored digraph whose uncolored digraph has $n+s$ vertices, $s \geq 0$ and consists of two cycles namely, one is an $n$-cycle and another one is an $(n-t)$-cycle, where $t \geq 1.$ For different values of $t$ as $t=1, 2$ and $3$ the bound on the $2$-exponent of the corresponding 2-colored digraphs, and the characterizations of the extremal 2-colored digraphs can be found in \cite{YubiiiGao:YyShao}, \cite{Shao:Gao:Sun} and \cite{Shao:Gao:Shanxi} respectively. In general, digraphs with an $n$-cycle and $(n-t)$-cycle ($t \geq 2$) have been considered in \cite{Huang:Liu}. In this paper the bound on the $2$-exponent and the corresponding extremal 2-colored digraphs have been characterized. Further, a much simpler bound on the $2$-exponent for this same family of digraphs has been obtained in  \cite{Shao:Gao:2}. In particular, for $s=0$ and $t \geq 3,$ authors in \cite{Shao:Gao:2} provided a tight upper bound on the $2$-exponent and characterized the extremal 2-colored digraphs. In 2006, S. Suwilo, B.L. Shader \cite{SuwWWw:Shaaadder} obtained an exact formula for $2$-exponent of $2$-primitive digraphs consisting of two arbitrary cycles. Using that formula they showed that the $2$-exponent of primitive extremal ministrong digraphs on $n$ vertices is at most $2n^2-7n+6.$ This is clearly an improvement of previously proposed bound of $2n^2-5n+3$ in \cite{SongGuu:Yaangg} where $2$-colored digraphs of the primitive ministrong digraphs with a given exponent are considered. In 2008 Y. Shao and Y. Gao \cite{Shao:Gao} considered $2$-colored digraphs of strongly connected symmetric digraphs with loops. In this paper they characterized $2$-primitive digraphs and found the sharp upper bound $3n-3$ on their exponents along with the characterization of the $2$-colored digraphs that attain the upper bound. Moreover, they established the exponent set $\{2,3, \ldots,3n-3\}$ for this family of digraphs on $n \geq 3$ vertices. Similar work has been done in \cite{Shao:Gao:1}, where the $2$-colored digraphs of strongly connected symmetric digraphs without loops are considered. Furthermore, S. Suwilo \cite{SaiSuWilo} proved that the $2$-exponent of an asymmetric 2-primitive $(n,s)$-lollipop on $n$ vertices with $s \leq n,$ is at most $\frac{(s^2-1)}{2}+(s+1)(n-s).$ The $(n,s)$-lollipops whose $2$-exponents achieving this bound are also characterized and for any asymmetric 2-primitive $(n,s)$-lollipop, a simple algorithm to find its 2-exponent can be found in this work.

In 1990 R.A. Brualdi and B. Liu. \cite{Bru:Liu} introduce the concept of three different kind of {\it generalized exponents} for a primitive digraph on $n$ vertices. In 2009, Y. Gao and Y. Shao \cite{YubinGaoYanlingShao} extend this concept to a primitive $k$-colored digraph and investigate the same for the primitive 2-colored Wielandt digraphs. In 2015, Mulyono, H. Sumardi, and S. Suwilo \cite{MulyonoSumardiSuwilo}, extend the notion of {\it scrambling index} of a primitive digraph to that of a 2-primitive digraph and discussed the same for 2-colored digraphs consisting of two cycles whose lengths differ by $1.$ The scrambling index for the class of 2-colored Hamiltonian digraphs on $n \geq 5$ odd vertices with two cycles of lengths $n$ and $\frac{(n-1)}{2}$ has been discussed in  \cite{Mardiningsih:Pasaribuuuuuu}. The {\it vertex exponent} or the {\it exponent of a vertex $v$}
in a 2-colored digraph $D^2$ is the smallest positive integer $s+t$ such that for each vertex $x$ in $D^2$ there is a walk of length $s+t$ consisting of $s$ red arcs and $t$ blue arcs. A few results on vertex exponent of 2-colored digraphs can be found in \cite{Saibsaibsaib}, \cite{AghniSaib}.

Similar work on digraphs with two cycles under various conditions, vertex exponent, generalized exponent, local exponent, (see \cite{Mardiningsih:Fathoni:Saib:1}), scrambling index for $2$-colored digraphs can be found in literature which are written in some other language than English. We expect that from the above discussion, reader will get a clear idea about almost all the works have been done in the context of $k$-primitivity.

{\bf Research problems:}
\begin{enumerate}
\item Let $k$ be a positive integer such that $k \in \bigg[\binom{n-1}{2}, \binom{n}{2}-\lceil \frac{n}{4} \rceil\bigg].$ Classify strongly connected $n$-tournaments which are $k$-primitive but not $(k+1)$-primitive.
\item Suppose $x_1, x_2, \ldots, x_l$ are positive integers. Consider the values of $x_1, x_2, \ldots, x_l$ except the values considered in \cite{Bea:Mou}. Does there exist an $(x_1, x_2, \ldots, x_l)$-game coloring of a strongly connected $n$-tournament $T$ that is $l$-primitive? What could be the result if we consider a general digraph $D$ instead of $T$?
\item Research can be focused on $k$-primitivity and $k$-exponent $(k \geq 2)$ for different kind of digraphs like digraphs with positive trace, regular graphs, Caylay digraphs, digraphs associated with doubly stochastic matrices etc. Such digraphs are not yet studied in the context of $k$-primitivity and results are not available in the literature. 
\end{enumerate}
\bibliographystyle{plain}
\bibliography{referexpo}

\begin{thebibliography}{10}

\bibitem{LeeRoyy:Kkirrk}
LeRoy~B. Beasley and Steve Kirkland.
\newblock A note on {$k$}-primitive directed graphs.
\newblock {\em Linear Algebra Appl.}, 373:67--74, 2003.
\newblock Special issue on the Combinatorial Matrix Theory Conference (Pohang,
  2002).

\bibitem{Bea:Mou}
LeRoy~B. Beasley and Sarah Mousley.
\newblock {$k$}-primitivity of digraphs.
\newblock {\em Linear Algebra Appl.}, 449:512--519, 2014.

\bibitem{Beassss:Neeealll}
LeRoy~B. Beasley and Cora~L. Neal.
\newblock Properties of 2-primitive tournament digraphs.
\newblock In {\em Proceedings of the {T}hirty-{F}ourth {S}outheastern
  {I}nternational {C}onference on {C}ombinatorics, {G}raph {T}heory and
  {C}omputing}, volume 163, pages 33--40, 2003.

\bibitem{Bru:Liu}
R.A. Brualdi and B.~Liu.
\newblock Generalized exponents of primitive directed graphs.
\newblock {\em J. Graph Theory}, 14(4):483--499, 1990.

\bibitem{Foorrn:Vaaaalll}
Ettore Fornasini and Maria~Elena Valcher.
\newblock Directed graph, {$2$}{D} state models, and characteristic polynomials
  of irreducible matrix pairs.
\newblock {\em Linear Algebra Appl.}, 263:275--310, 1997.

\bibitem{Forn:Val}
Ettore Fornasini and Maria~Elena Valcher.
\newblock Primitivity of positive matrix pairs: algebraic characterization,
  graph-theoretic description, and 2{D} systems interpretation.
\newblock {\em SIAM J. Matrix Anal. Appl.}, 19(1):71--88, 1998.

\bibitem{Frob:2}
G.~Frobenius.
\newblock \"{U}ber matrizen aus nicht negativen elementen.
\newblock {\em S. -B. K. Preuss. Akad. Wiss. Berlin}, pages 456--477, 1912.

\bibitem{YubiiiGao:YyShao}
Yubin Gao and Yanling Shao.
\newblock Exponents of two-colored digraphs with two cycles.
\newblock {\em Linear Algebra Appl.}, 407:263--276, 2005.

\bibitem{YubinGaoYanlingShao}
Yubin Gao and Yanling Shao.
\newblock Generalized exponents of primitive two-colored digraphs.
\newblock {\em Linear Algebra Appl.}, 430(5-6):1550--1565, 2009.

\bibitem{Gaooo:Shaaaoooo}
Yubin Gao and Yanling Shao.
\newblock On {$k$}-primitivity of two classes of digraphs.
\newblock {\em Linear Algebra Appl.}, 430(11-12):2922--2928, 2009.

\bibitem{Heap:Lynn}
B.~R. Heap and M.~S. Lynn.
\newblock The index of primitivity of a non-negative matrix.
\newblock {\em Numer. Math.}, 6:120--141, 1964.

\bibitem{Herstein}
I.~N. Herstein.
\newblock A note on primitive matrices.
\newblock {\em Amer. Math. Monthly}, 61:18--20, 1954.

\bibitem{Huang:Liu}
Fengying Huang and Bolian Liu.
\newblock Exponents of a class of two-colored digraphs with two cycles.
\newblock {\em Linear Algebra Appl.}, 429(2-3):658--672, 2008.

\bibitem{SongGuu:Yaangg}
Sang~Gu Lee and Jeong~Mo Yang.
\newblock Bound for 2-exponents of primitive extremal ministrong digraphs.
\newblock {\em Commun. Korean Math. Soc.}, 20(1):51--62, 2005.

\bibitem{Liu:Lai}
Bolian Liu and Hong-Jian Lai.
\newblock {\em Matrices in combinatorics and graph theory}, volume~3 of {\em
  Network Theory and Applications}.
\newblock Kluwer Academic Publishers, Dordrecht, 2000.
\newblock With a foreword by Richard A. Brualdi.

\bibitem{Mardiningsih:Fathoni:Saib:1}
Mardiningsih, Muhammad Fathoni, and Saib Suwilo.
\newblock Local exponents of two-colored bi-cycles whose lengths differ by $1$.
\newblock In {\em Special Issue: The 10th IMT-GT International Conference on
  Mathematics, Statistics and its Applications 2014}, Malaysian Journal of
  Mathematical Sciences, pages 205--218, 2016.

\bibitem{Mardiningsih:Pasaribuuuuuu}
Mardiningsih and Merryanty~L. Pasaribu.
\newblock The scrambling index of a class of two-colored hamiltonian digraphs.
\newblock In {\em International Conference on Mathematics, Engineering and
  Industrial Applications, AIP Conf. Proc. 1775, (2016)}, pages
  030056--1--030056--8, 2016.

\bibitem{Minc}
Henryk Minc.
\newblock {\em Nonnegative matrices}.
\newblock Wiley-Interscience Series in Discrete Mathematics and Optimization.
  John Wiley \& Sons, Inc., New York, 1988.
\newblock A Wiley-Interscience Publication.

\bibitem{MulyonoSumardiSuwilo}
Mulyono, H.~Sumardi, and Saib Suwilo.
\newblock The scrambling index of primitive two-colored two cycles whose
  lengths differ by $1$.
\newblock {\em Far East J. Math. Sci.}, 96(1):113--132, 2015.

\bibitem{Ole:Shad:Dri}
D.~D. Olesky, Bryan Shader, and P.~van~den Driessche.
\newblock Exponents of tuples of nonnegative matrices.
\newblock {\em Linear Algebra Appl.}, 356:123--134, 2002.
\newblock Special issue on algebraic graph theory (Edinburgh, 2001).

\bibitem{Per}
Peter Perkins.
\newblock A theorem on regular matrices.
\newblock {\em Pacific J. Math.}, 11:1529--1533, 1961.

\bibitem{Perr}
Oskar Perron.
\newblock Zur {T}heorie der {M}atrices.
\newblock {\em Math. Ann.}, 64(2):248--263, 1907.

\bibitem{Port}
V.~Yu. Protasov.
\newblock Classification of {$k$}-primitive sets of matrices.
\newblock {\em SIAM J. Matrix Anal. Appl.}, 34(3):1174--1188, 2013.

\bibitem{Prot:Voy}
V.~Yu. Protasov and A.~S. Voynov.
\newblock Sets of nonnegative matrices without positive products.
\newblock {\em Linear Algebra Appl.}, 437(3):749--765, 2012.

\bibitem{Romanovskyyy}
V.~Romanovsky.
\newblock Un th\'eor\`eme sur les z\'eros des matrices non n\'egatives.
\newblock {\em Bull. Soc. Math. France}, 61:213--219, 1933.

\bibitem{Shad:Suw}
Bryan~L. Shader and Saib Suwilo.
\newblock Exponents of nonnegative matrix pairs.
\newblock {\em Linear Algebra Appl.}, 363:275--293, 2003.
\newblock Special issue on nonnegative matrices, $M$-matrices and their
  generalizations (Oberwolfach, 2000).

\bibitem{Shao:Gao}
Yanling Shao and Yubin Gao.
\newblock Exponents of 2-coloring of symmetric digraphs.
\newblock {\em Linear Algebra Appl.}, 428(7):1538--1550, 2008.

\bibitem{Shao:Gao:1}
Yanling Shao and Yubin Gao.
\newblock Exponents of 2-colorings of loopless, symmetric digraphs.
\newblock {\em Linear Multilinear Algebra}, 57(1):65--74, 2009.

\bibitem{Shao:Gao:Shanxi}
Yanling Shao and Yubin Gao.
\newblock Exponents of two-colored digraphs.
\newblock {\em Czechoslovak Math. J.}, 59(134)(3):655--685, 2009.

\bibitem{Shao:Gao:2}
Yanling Shao and Yubin Gao.
\newblock On the exponents of two-colored digraphs with two cycles.
\newblock {\em Linear Multilinear Algebra}, 57(2):185--199, 2009.

\bibitem{Shao:Gao:Sun}
Yanling Shao, Yubin Gao, and Liang Sun.
\newblock Exponents of a class of two-colored digraphs.
\newblock {\em Linear Multilinear Algebra}, 53(3):175--188, 2005.

\bibitem{Shen:9}
J.~Shen.
\newblock {\em Exponents of primitive digraphs}.
\newblock 1998.
\newblock Thesis (Ph.D.)--Queen's University, Kingston, Ontario, Canada.

\bibitem{Suwiloooos}
Saib Suwilo.
\newblock {\em On 2-exponents of 2-digraphs}.
\newblock ProQuest LLC, Ann Arbor, MI, 2001.
\newblock Thesis (Ph.D.)--University of Wyoming.

\bibitem{SaiSuWilo}
Saib Suwilo.
\newblock $2$-exponents of two-coloured lollipops.
\newblock {\em MATEMATIKA}, 21(1):11--22, 2008.

\bibitem{Saibsaibsaib}
Saib Suwilo.
\newblock vertex exponents of two-colored primitive extremal ministrong
  digraph.
\newblock {\em Global journal of technology and optimization}, 2:166--174,
  2011.

\bibitem{SuwWWw:Shaaadder}
Saib Suwilo and Bryan~L. Shader.
\newblock On $2$-exponents of ministrong $2$-digraphs.
\newblock In {\em Proceedings of the $2^{\begin{tiny} nd \end{tiny}}$ IMT-GT
  Regional Conference on Mathematics, Statistics and Applications}, pages
  51--56, 2006.

\bibitem{AghniSaib}
Aghni Syahmarani and Saib Suwilo.
\newblock Vertex exponents of a class of two-colored {H}amiltonian digraphs.
\newblock {\em J. Indones. Math. Soc.}, 18(1):1--19, 2012.

\bibitem{Thhhom}
G.L. Thompson.
\newblock {\em Lectures on Game Theory, Markov Chains, and Related Topics}.
\newblock Sandia Corporation Monograph SCR-11, 1958.

\bibitem{Wie}
Helmut Wielandt.
\newblock Unzerlegbare, nicht negative {M}atrizen.
\newblock {\em Math. Z.}, 52:642--648, 1950.

\end{thebibliography}
\end{document}